\documentclass{amsart}
\usepackage{amsmath}
\usepackage{amssymb}
\newtheorem{Proposition}{Proposition}
  \newtheorem{Remark}{Remark}
  \newtheorem{Corollary}[Proposition]{Corollary}
  \newtheorem{Lemma}[Proposition]{Lemma}
   \newtheorem{Theorem}[Proposition]{Theorem}

\newtheorem{Note}[Proposition]{Note}

\def\z{\noindent}
\def\bchi{\ensuremath{\mbox{\large $\chi$}}}
\def\lap{\mathcal{L}}
\def\bor{\mathcal{B}}
\def\CC{\mathbb{C}}

\def\RR{\mathbb{R}}

\def\NN{\mathbb{N}}
\def\erm{\mathrm{e}}

\begin{document}

\title{Nonlinear Stokes phenomena in first or second order differential equations\\
  ~~\\~~\\{\sc Dedicated to Professor Kawai on the occasion of his 60th
    birthday}} \author{O.  Costin}
\address{Mathematics Department \\The Ohio State University\\
  Columbus OH 43210\\USA\\ costin@math.ohio-state.edu}

\date{}
\bigskip

\begin{abstract} We study singularity formation in nonlinear differential equations of order $m\leqslant 2$,
  $y^{(m)}=A(x^{-1},y)$. We assume $A$ is analytic at $(0,0)$ and $\partial_y
  A(0,0)=\lambda\ne 0$ (say, $\lambda=(-1)^m$). If $m=1$ we assume
  $A(0,\cdot)$ is meromorphic and nonlinear. If $m=2$, we assume $A(0,\cdot)$
  is analytic except for isolated singularities, and also that
  $\int_{s_0}^\infty |\Phi(s)|^{-1/2}d|s|<\infty$ along some path avoiding the
  zeros and singularities of $\Phi$, where $\Phi(s)=\int_{0}^s A(0,\tau)d\tau$.
  Let $H_{\alpha}=\{z:|z|>a>0,\arg(z)\in (-\alpha,\alpha)\}$.
  
  If the Stokes constant $S^+$ associated to $\RR^+$ is nonzero, we show that
  {\em all} $y$ such that $\lim_{x\to +\infty}y(x)=0$ are singular at $2\pi
  i$-quasiperiodic arrays of points near $i\RR^+$. The array location
  determines and is determined by $S^+$.  Such settings include the Painlev\'e
  equations $P_I$ and $P_{II}$.  If $S^+=0$, then there is exactly one
  solution $y_0$ without singularities in $H_{2\pi-\epsilon}$, and $y_0$ is
  entire iff $y_0=A(z,0)\equiv 0$.
  
  The singularities of $y(x)$ mirror the singularities of the Borel transform
  of its asymptotic expansion, $\mathcal{B}\tilde{y}$, a nonlinear analog of
  Stokes phenomena.  If $m=1$ and $A$ is a nonlinear polynomial with
  $A(z,0)\not\equiv 0$ a similar conclusion holds even if $A(0,\cdot)$
  is linear. This follows from the property that if $f$ is superexponentially
  small along $\RR^+$ and analytic in $H_{\pi}$, then $f$ is
  superexponentially unbounded in $H_{\pi}$, a consequence of decay estimates
  of Laplace transforms.

  Compared to \cite{Inventiones} this analysis is restricted to first and
  second order equations but shows that singularities {\em always} occur, and
  their type is calculated in the polynomial case.  Connection to
  integrability and the Painlev\'e property are  discussed.

\end{abstract}
\gdef\shorttitle{Nonlinear Stokes phenomena}
\maketitle
\setcounter{section}{0}
\section{Introduction} 

Stokes transitions play an important role in the study of differential
equations (see \cite{Kawai}, \cite{Its}, \cite{Sibuya} and references
therein). In linear differential equations, solutions that are small in some
sector usually become exponentially large in complementary ones. 

There is an interesting analog for analytic nonlinear equations which we study
for first and second order ones.

Assume that $A(z,y)$ is analytic at $(0,0)$, and that the following
nondegeneracy condition holds: $\partial_y A(0,y):=\lambda\ne 0$. Let, for
$m=1,2$,
\begin{equation}
    \label{eq2}
   y^{(m)}=A(1/x,y) 
  \end{equation}
  By a simple change of variable $y\mapsto y+a+b x^{-1}+cx^{-2}$ we can
  assume, without loss of generality, $A(z,0)=o(z^2)$.  The conditions in
  \cite{DMJ} and \cite{Inventiones} (which do not require $m\leqslant 2$) then
  apply. By a change of independent variable we can make $\lambda=(-1)^m$.
  Let $A_1'(0)=\alpha $.
  
  There exists a one parameter family of solutions of (\ref{eq2}) which decay
  algebraically as $x\to\infty$ \cite{Wasow}, \cite{DMJ}. All these solutions
  can be written as  generalized Borel summed transseries \cite{DMJ} (the
  notation here is similar to that in \cite{Inventiones})
\begin{equation}
    \label{eq3}
 {y}(x)=
\sum_{{k}\ge 0}
  {C}^{{k}}\erm^{-{k}x}
  x^{\alpha k }{y}_{{k}}(x)\\ 
=\sum_{{k}\ge 0}
  {C}^{{k}}\erm^{-{k}x}
  x^{\alpha k}\mathcal{L}\mathcal{B}\tilde{{y}}_{{k}}(x)\equiv\mathcal{L}\mathcal{B}\tilde{{y}}(x)
\end{equation}
where $\tilde{{y}}_{k}$ are formal series in integer powers of $x^{-1}$. The
functions $y_k:=\lap\bor\tilde{y}_k$, the generalized Borel sums of
$\tilde{y}_k$, are analytic for $\Re(x)>\nu>0$ uniformly in $k$. The function
series in (\ref{eq3}) is absolutely and uniformly convergent for $\Re(x)>\nu$
and defines an analytic function there.  Any decaying solution can be
represented uniquely in terms of Laplace transforms in the upper half plane,
$y_k^+=\lap(\bor \tilde{y})^+$ where $^+$ indicates the upper half plane
continuation or in terms of the lower half plane continuations \cite{DMJ}
\begin{equation}
    \label{eq31}
 {y}(x)=\sum_{{k}\ge 0}
{C^+}^k e^{-{k}x}
  x^{\alpha k}y_k^+=\sum_{{k}\ge 0}{C^-}^ke^{-{k}x}
  x^{\alpha k}y_k^-
\end{equation}
and we have (\cite{DMJ} Eq. (1.27))
\begin{equation}
  \label{eq:c+-}
  C^+-C^-=S^+
\end{equation}
where $S^+=S_1$ in \cite{DMJ} is the Stokes constant associated to $\RR^+$.

In this paper we prove, under some further assumptions, that all solutions,
except at most one, that go to zero along some direction form arrays of
singularities near $i\RR$.  After normalization of the equation, the array
location is given by
\begin{equation}
    \label{eq702}
x_n=2n\pi i+\alpha\ln(2n\pi i)+\ln C^{\pm}-\ln\xi_s+o(1) \ \ (\NN\ni n\to\pm \infty)
\end{equation} for some $\xi_s\ne 0$. (If
$C^-=0$ ($C^+=0$), then $y$ has no singularities in the lower (upper,
respectively) half plane.)  In \cite{Inventiones}, Eq. (28) we showed that
such arrays occur provided a solution of an associated equation, (\ref{eq5})
below, is singular. Here we show that the solution is singular indeed, under
our assumptions.

There is a correspondence between singularities in inverse Laplace space and
singularities in the original space, a nonlinear analog of Stokes phenomena:
if we change variables to make $\partial_y A(0,0)=\lambda$, then $n$ becomes
$n/\lambda$. Singularities in Borel space are spaced by (exactly) $|\lambda|$,
see \cite{DMJ}, while in the $x$ plane their mirror singularities are, up to
log corrections, as we see from (\ref{eq702}), regularly spaced by
$|\lambda|^{-1}$.

Part of the results rely on a general property that if an analytic function is
superexponentially small along a ray, then it becomes superexponentially large
along a complementary ray.

\section{Main results}\label{Singf}
\z In \cite{Inventiones} it is shown that there is a unique analytic solution
$F_0$ of the associated equation
\begin{equation}
    \label{eq5}
\sum_{k=1}^m \xi^k F_0^{(k)}=A(0,F_0)
\end{equation}
with $\xi^{-1} F_0(\xi)\to 1$ as $\xi\to 0$. $F_0$ has the property that
\begin{equation}
  \label{eq:asympty}
 y(x)=F_0(x^{\alpha} e^{-x})+O(x^{-1})
\end{equation}
as $x\to\infty$ if $\xi=x^{\alpha} e^{-x}$ stays bounded, and if $\xi_s$ is a
singular point of $F_0$, then $y$ is singular along the array (\ref{eq702}).
From (\ref{eq:asympty}) and the argument principle, if $F_0(\xi)$ has a pole
or an algebraic or logarithmic branch point at $\xi_s$ then the actual
solution $y$ of (\ref{eq2}) has a quasiperiodic array of singularities, to
leading order of the same type as that of $F_0$.  The origin of the
singularities of $y(x)$ is the fact that in (\ref{eq3}) the exponentials
$e^{-kx}$ become $O(1)$ in the antistokes direction, and they can ``pile up''
to create blow-up of the solution.  The singularities form $2\pi
i$-quasiperiodic arrays since they are governed by the size of
$x^{\alpha}e^{-x}$. If $y^-$ is the solution analytic in the lower half plane,
then $C^-=0$ and $S^+$ can be calculated from (\ref{eq702}) and
(\ref{eq:c+-}).

\z {\bf Assumptions}. We recall $A$ is analytic at $(0,0)$. Let
$H=A(0,\cdot)$ and $\Phi(s)=\int_{0}^s H(\tau)d\tau$. If $m=1$ assume $H$ is
meromorphic and nonlinear. If $m=2$ assume $H$ is analytic except for isolated
singularities and such that $\int_a^\infty |\Phi(s)|^{-1/2}ds<\infty$ along a
path $P$ avoiding the zeros and singularities of $\Phi$. This includes the
Painlev\'e equations $P_I$ and $P_{II}$, where $H$ is a nonlinear polynomial
\cite{Inventiones}.  Consider the one-parameter family $\mathcal{F}$ of
solutions $y$ of (\ref{eq2}) such that $\lim_{x\to +\infty}y(x)=0$.
\begin{Theorem}\label{Lexp1}
  (i) Under the assumptions above, all functions in $\mathcal{F}$ except at
  most one are singular in the upper or lower half--plane.  The singularities
  form arrays given by (\ref{eq702}).  If $H$ is polynomial of degree $N>1$,
  then the local behavior of $y$ near a singular point $x_s$ is to leading
  order $(x-x_s)^{m/(1-N)}$.
  
  (ii) If there exists an $y_0\in\mathcal{F}$ analytic in a sector of opening
  more than $\pi$ centered on $\RR^+$, then there is no other such $y_0$, and
   $S^+=0$; $y_0$ is entire iff $y_0\equiv A(x,0)\equiv 0$.
  
   (iii) For $m=1$, if $A(z,y)$ is polynomial of degree $N>1$ in $y$ and
   $A(z,0)\not\equiv 0$, all solutions in $\mathcal{F}$ except at most one
   have infinitely many singularities in $H_{\pi+\epsilon}$, even if
   $A(0,\cdot)$ is linear.  The singularity type is as in (i).
\end{Theorem}

\subsection{General decay estimates} 
We state separately some useful and relatively basic estimates, helpful in
proving Theorem~\ref{Lexp1} (iii).\footnote{ Lemma~\ref{Carson+} resembles
  Carlson's Lemma \cite{T} pp. 185, but does not appear to imply it or to
  follow from it. Also, the referee points out that similar estimating methods
  appear in \cite{PW}.}
\begin{Lemma}\label{Carson+}
  If the function $f\not\equiv 0$ is analytic in a half plane
  $\{z:\Re(z)>\beta\}$ and for all $a>0$ we have $f(t)=O(e^{-at})$ for
  $t\in\RR^+$, $t\to +\infty$, then for any $\epsilon \geqslant 0$ small
  enough $ f(z) e^{-bz^{1-\epsilon}}$ is unbounded in the closed right half
  plane for all $b>0$.
\end{Lemma}

\begin{Lemma}\label{Lexp}
  Assume $F\in L^1(\RR^+)$ and for some $\epsilon>0$ we have
  \begin{equation}
    \label{eq:estlower1}
     \mathcal{L}F(x)=O(e^{-\epsilon x})\ \ \ \text{as}\ \ \ x\to +\infty
  \end{equation}
  Then $F=0$ a.e. on $[0,\epsilon]$. (The result is sharp as discussed after
  the proof.)
\end{Lemma}
\begin{Corollary}\label{decay11}
  Assume $F\in L^1$ and $\mathcal{L}F=O(e^{-ax})$ as $x\to +\infty$ for all
  $a>0$. Then $F=0$ a.e. on $\RR^+$.
\end{Corollary}
\begin{Corollary}\label{Ast}
  Assume $\alpha>\frac{1}{2},f\not\equiv 0$ is analytic in $S=\{z: |z| >R,
  2|\arg z|\leqslant \pi/\alpha\}$ and for any $a>0$ we have $f(t)=o(e^{- a
    t^\alpha})$ as $t\to +\infty$. Then for any $b,p>0$ and
  $0<\alpha'<\alpha$, the function $f(z)e^{-{ bz^{{\alpha}'}}}$ is unbounded
  in $S$.
\end{Corollary}

\section{Proofs}
\subsection{Proof of Theorem~\ref{Lexp1}} (i) Let first $m=1$. We show by
contradiction that $F_0$ in (\ref{eq5}) is not entire. Assume $g(t):=F_0(e^t)$
was entire. We have $g'=H(g)$.  First, $H$ must have roots or else $G=1/H$
would be entire and then the function $Q$ given by
\begin{equation}
  \label{eq:eqq2}
  Q(g)=\int{G(g)dg}+C=t
\end{equation}
would be entire, and injective since $Q$ has an (entire) inverse. Thus $Q$ is
linear.

Assume, without loss of generality, that $H(0)=0$. Next we show that
if $H(\tau)=0$ then $g(t_0)\ne\tau$ for all $t_0\in\CC$.  Indeed, if not, $g'(t_0)=0$ and we would
have
\begin{equation}
  \label{eq:eqq}
  \frac{(g(t)-\tau)'}{g(t)-\tau}=\frac{H(g(t))-H(\tau)}{g(t)-\tau}\to
  H'(\tau)\ \ \text{as} \ t\to t_0
\end{equation}
But this is a contradiction, since $t_0$ is a zero of finite order of
$g(t)-\tau$ and the left side of (\ref{eq:eqq}) tends to infinity.  Since
$g(t)\ne 0$ for all $t\in\CC$, we have $g=e^h$ with $h$ entire, and since $g$
is nontrivial it cannot avoid any other value, thus $0$ is the only root of
$H$.  Then $h'=e^{-h}H(e^h)$ where now the right hand side is meromorphic and
everywhere nonzero, a case we have already analyzed.

Let now $m=2$, $g(t)=F_0(e^{-t})$. We have $g''=H(g)$.  Then, multiplying by
$g'$ and integrating once, we get ${g'}^2=\Phi(g)+C$. Using the condition
$g\sim e^{-t}$ for large $t$ and the analyticity of $H$ for small $g$ we get
$C=0$.  Choosing a determination of the log and a $t_0$ such that $g(t_0)$ is
not a zero or a singularity of $\Phi$ we get
\begin{equation}
  \label{eq:eq72}
J(g(t_0),g(t)):=\int_{g(t_0)}^{g(t)}\Phi(s)^{-1/2}ds=t-t_0
\end{equation}
We take the path $P$ towards infinity. We have that
$J=\lim_{g\to\infty}J(g(t_0),g)$, $g\in P$, is finite and $t_s:=t_0+J$
is clearly a singular point of $g$.

(ii) It follows from (\ref{eq:c+-}) that unless $S^+=0$ every solution
develops singularities in the upper or lower half plane (in both, except when
$C_-=0$ or $C^+=0$). Conversely, if $S^+=0$, it follows from \cite{DMJ} that
$\bor\tilde{y}$ is analytic in the $\CC\setminus\RR^-$ (Theorem 1) and
exponentially bounded in distributions in $H_{2\pi-\epsilon}$ (Theorem 2).
Thus $\lap\bor\tilde{y}$ is a solution of the original differential equation
with the asymptotic expansion $\tilde{y}$ in any sector of opening less than
$2\pi+\pi/2$ centered on $\RR^+$. In particular, this solution cannot be
entire unless it is identically zero, and then $A_0\equiv 0$.  The type of
singularity follows from (\ref{eq:eq72}) and the discussion at the beginning
of \S\ref{Singf}.

For part (iii) we need the following result.
 \begin{Lemma}\label{Formsing}
    Assume $y$ solves
\begin{equation}
    \label{eq22}
y'=\sum_{k=0}^N B_k(x)y^k
\end{equation}
$N>1$, in a sector $S$ where the analytic coefficients $B_k$ are $O(x^{q_k})$
and $B_N\not\equiv 0$.  Assume furthermore that for a sequence $x_j\to\infty$
in $S$, $y(x_j) $ grows faster than polynomially in $x_j$. Then for all $j$
large enough, within a distance $O(x_j^{-q_N} y(x_j)^{1-N})$ of $x_j$ there is
a singularity of $y$.
  \end{Lemma}
  \z {\bf Note.} Heuristically, near a large $x_j$, the dominant balance is of
  the form $y'\sim C y^N x^{q_N},\, C\ne 0$, which forms a singularity as described.
  \begin{proof}Let $x=x_j+\sigma$ and $\eta=1/y$. The equation for $\sigma$
    is
\begin{equation}
  \label{eq:eqsigma}
  \frac{d\sigma}{d\eta}=-\frac{\eta^{N-2}}{\sum_{k=0}^N\eta^k B_{N-k}(x_j+\sigma)}=G(\sigma,\eta)
\end{equation}
with the initial condition $\sigma(1/y(x_j))=0$. It is straightforward to
check that in our assumptions, for large enough $x_j$, $G(\sigma,\eta)$ is
analytic in the polydisk $\{(\sigma,\eta):|\sigma|<1,|\eta|<2/|y(x_j)|\}$
where furthermore $|G(\sigma,\eta)|=O(\eta^{N-2}x_j^{-q_N})$. Then, in the disk
$\{\eta:|\eta|< 2/|y(x_j)|\}$ there exists a unique analytic solution
$\sigma$ and $\sigma=O(\eta^{N-1}x_j^{-q_N})$. It is clear that $x_j+\sigma(0)$
is a zero of $\eta$ and thus a singularity of $y$.
\end{proof}
\z {\em Proof of Theorem~\ref{Lexp1} (iii)} Assume that there is a solution
$y_0$ of such that $y_0\to 0$ as $x\to+\infty$ with finitely many
singularities in a sector $S$ of opening $\pi+2\epsilon$ centered on $\RR^+$.
(There may be an exceptional solution, for instance when the Stokes constant
is zero and the asymptotic series converges.)  Let $y_2$ be another solution
of $y'=A(x^{-1},y)$ such that $y_2\to 0$ as $x\to\infty$ and let
$\delta=y_2-y_1$.  The equation for $\delta$ is of the form (\ref{eq22}) with
$B_0=0$. From the general theory of differential equations \cite{Wasow} or
from \cite{DMJ} it follows that $\delta=O(e^{-x})$ as $x\to\infty$. By
Lemma~\ref{Carson+} there is a sequence of $x_j\to \infty$ in $S$ so that
$\delta\ge const.  \exp(|x_j|^{1-\epsilon})$.  Near every $x_j$ with $j$ large
enough there is, by Lemma~\ref{Formsing}, a singularity of $\delta$ and thus
of $y_2$.
\subsection{Connection with integrability} It is seen that if $H$ is a
polynomial of degree $N>2$ if  $m=1$ or $N> 3$ if $m=2$ a one
parameter family of solutions forms arrays of branch-point singularities in
$\CC$ and the Painlev\'e test fails. For other values of $p$, the corrections
$F_j,j>1$ \cite{Inventiones} can be calculated to determine the exact type of
singularity. More can be done however. Assuming we are in the case where
singularities are branch points, in view of the exact description of type and
location of singularities it is possible to calculate the monodromy of
solutions along a curve in $\CC$ winding among sufficiently many singularities
to show dense branching, a concept proposed by Kruskal, which can be
used to show absence of continuous first integrals as in \cite{Rodica1},
\cite{Rodica2}. This will be the subject of a different paper.
\subsection{Proof of Lemma~\ref{Lexp}}\label{L+}
  We write
  \begin{equation}
    \label{eq:splitlap}
    \int_0^\infty e^{-px}F(p)dp=\int_0^\epsilon e^{-px}F(p)dp+\int_\epsilon^\infty e^{-px}F(p)dp
  \end{equation}
we note that
\begin{equation}
    \label{eq:splitlap2}
    \Big|\int_\epsilon^\infty e^{-px}F(p)dp\Big|\leqslant 
     e^{-\epsilon x}\int_\epsilon^\infty |F(p)|dp\leqslant  e^{-p\epsilon}\|F\|_1=O(e^{-\epsilon x})
  \end{equation}
Therefore
\begin{equation}
    \label{eq:estlower3}
    g(x)=\int_0^\epsilon e^{-px}F(p)dp =O(e^{-\epsilon x})\ \ \ \text{as}\ \ \ x\to +\infty
  \end{equation}
  The function $g$ is manifestly entire. Let $h(x)=e^{\epsilon x}g(x)$. Then
  by assumption $h$ is entire and uniformly bounded for $x\in\RR$ (since by
  assumption, for some $x_0$ and all $x>x_0$ we have $|h|\leqslant C$ and by
  continuity $\max|h|<\infty$ on $[0,x_0]$).  The function is bounded by
  $\|F\|_1$ for $x\in i\RR$.  By Phragm\'en-Lindel\"of's theorem (first
  applied in the first quadrant and then in the fourth quadrant, with
  $\beta=1, \alpha=2$) $h$ is bounded in the closed right half plane. Now, for
  $x=-s<0$ we have
\begin{equation}
  \label{eq:estlower4}
  e^{-s\epsilon}\int_0^\epsilon e^{sp} F(p) dp\leqslant  \int_0^\epsilon | F(p)|\leqslant  \|F\|_1
\end{equation}
Again by Phragm\'en-Lindel\"of (and again applied twice) $h$ is bounded in the
closed left half plane thus bounded in $C$, and it is therefore a constant.
But, by the Riemann-Lebesgue lemma, $h\to 0$ for $x=is$ when $s\to +\infty$.
Thus $h\equiv 0$. Therefore, with $\bchi_A$ the characteristic function of $A$,
\begin{equation}
  \label{eq:fourier2}
  \int_0^\epsilon F(p) e^{-isp}dp=\hat{\mathcal{F}}(\bchi_{[0,\epsilon]}F)=0
\end{equation}
for all $s\in\RR$ entailing the conclusion.

\z {\bf Note}. 
In the opposite direction, by Laplace's method it is easy to check
that for any small $\epsilon>0$ we have $\displaystyle \mathcal{L}
e^{-p^{-\frac{1-\epsilon}{\epsilon}}} =o\left(e^{-{x^{1-\epsilon}}}\right)$
and for any $n$ $\displaystyle \mathcal{L} \left(e^{- E_{n+1}(1/p)}\right)=o(e^{-{x/L_n(x)}})$
where $E_n$ is the n-th composition of the exponential with itself and $L_n$
is the n-th composition of the log with itself.
$\Box$.
\subsection{Proof of Proposition~\ref{Carson+}}\label{C+}
By a change of variable we may assume that $\beta=-1$. Assume that for some
$b>0$ $e^{-bz^{1-\epsilon}}f(z)$ was bounded in the closed right half plane.
Then $ \psi(z)= (1+z)^{-2}{e^{-bz^{1-\epsilon}}f(z)} $ satisfies the
assumptions of Lemma~\ref{Invlap}. But then $\psi(z)= \mathcal{L}
\mathcal{L}^{-1}\psi(z)$ satisfies the assumptions of Corollary \ref{decay11}
and $\psi\equiv 0$.
\begin{Note}
  There is indeed loss of exponential rate: the entire function
  $\Gamma(x)^{-1}=O(e^{-n\ln n})$ on $\RR^+$ but is bounded by
  $O(e^{|x|\pi/2})$ in the closed right half plane.
\end{Note}
\section{Appendix: overview of Laplace transform properties}
\z For convenience we  provide some standard results on Laplace transforms.
\begin{Remark}[Uniqueness] \label{uniqLap}
  Assume $F\in L^1(\RR^+)$ and $\mathcal{L}F=0$ for a set of $x$  with an
  accumulation point. Then $F=0\,\, a.e.$
\end{Remark}
\begin{proof}
   By analyticity, $\mathcal{L}F=0$ in the open right half plane and by
continuity, for $s\in\RR$, $\mathcal{L}F(is)=0=\hat{\mathcal{F}} F$ where
$\hat{\mathcal{F}} F$ is the Fourier transform of $F$ (extended by zero for
negative values of $p$).  Since $F\in L^1$ and $0=\hat{\mathcal{F}} F\in L^1$,
by the known Fourier inversion formula \cite{Rudin}, $F=0$.
\end{proof}
\begin{Lemma}\label{Invlap}
  Assume that $ c\geqslant  0$ and  $f(z)$ is analytic in 
  $H_c:=\{z:\Re \, z\geqslant  c\}$.  Assume further that  $g(t):=\sup_{c'\geqslant 
    c}|f(c'+it)|\in L^1(\RR,dt)$.  Let
\begin{equation}
  \label{eq:invlp1}
  F(p)=\frac{1}{2\pi i}\int_{c-i\infty}^{c+i\infty}e^{px}f(x)dx=
:(\mathcal{L}^{-1}f)(p)
\end{equation}
Then for any $x\in \{z:\Re \, z> c\}$ we have $\displaystyle
\mathcal{L}F=\int_0^{\infty}e^{-px}F(p)dp=f(x)$
\end{Lemma}
Note that for any $x'=x'_1+iy'_1\in \{z:\Re \, z> c\}$
\begin{equation}
  \label{eq:estFub}
  \int_0^{\infty}dp\int_{c-i\infty}^{c+i\infty} \left|e^{p(s-x')}f(s)\right|d|s|
\leqslant  \int_0^{\infty}dp e^{p(c-x'_1)}\|g\|_1\leqslant 
\frac{\|g\|_1}{x'_1-c}
\end{equation} 
and thus, by Fubini we can interchange the orders
  of integration:
  \begin{multline}
    \label{eq:eqU}
    U(x')=\int_0^{\infty}e^{-px'}\frac{1}{2\pi
    i}\int_{c-i\infty}^{c+i\infty}e^{px}f(x)dx\\=\frac{1}{2\pi
    i}\int_{c-i\infty}^{c+i\infty}dx f(x)\int_0^{\infty}dp e^{-px'+px}=
\frac{1}{2\pi
    i}\int_{c-i\infty}^{c+i\infty} \frac{f(x)}{x'-x}dx
  \end{multline}
  Since $g\in L^1$ there exist subsequences $\{\tau_n\},\{-\tau'_n\}$ tending
  to infinity 
  such that $|g(\tau_n)|\to 0$. Let $x'>\Re x=x_1$ and consider the box
  $B_n=\{z:\Re z \in[x_1,x'], \Im z\in [-\tau'_n,\tau_n]\}$ with positive
  orientation.

  \begin{equation}
  \label{eq:Bn}
  \int_{B_n}\frac{f(s)}{x'-s}ds=-f(x')
\end{equation}
while, by construction,
\begin{equation}
  \label{eq:Bn1}
 \lim_{n\to\infty} \int_{B_n}\frac{f(s)}{x'-s}ds=\int_{x'-i\infty}^{x'+i\infty} \frac{f(s)}{x'-s}ds-\int_{c-i\infty}^{c+i\infty} \frac{f(s)}{x'-s}dx
\end{equation}
On the other hand, by dominated convergence, we have
\begin{equation}
  \label{eq:Bn3}
  \int_{x'-i\infty}^{x'+i\infty} \frac{f(s)}{x'-s}ds\to 0\ \ \text{as}\ \ x'\to\infty
\end{equation}

\z {\bf Acknowledgments}. The work was partially
supported by NSF grant 0406193.

\end{document}